\documentclass{birkjour}
\usepackage{xcolor, hyperref, bbold}
 \newtheorem{thm}{Theorem}[section]
 
 \newtheorem{lem}[thm]{Lemma}
 \newtheorem{prop}[thm]{Proposition}
 \theoremstyle{definition}
 \newtheorem{defn}[thm]{Definition}
 \theoremstyle{remark}
 \newtheorem{rem}[thm]{Remark}
 \newtheorem*{ex}{Example}
 \numberwithin{equation}{section}

\begin{document}

%
%
%
%
%
%
%
%
%

\title[Brownian motion, martingales and It\^o formula in Clifford analysis]
 {Brownian motion, martingales and It\^o formula in Clifford analysis}

\author[S. Bernstein]{Swanhild Bernstein}

\address{%
Institute for Applied Analysis\\
TU Bergakademie Freiberg\\
09599 Freiberg \\
Germany}
\email{Swanhild.Bernstein@math.tu-freiberg.de}

\author[D. Legatiuk]{Dmitrii Legatiuk}
\address{Chair of Mathematics\\
Universit\"at Erfurt \\
99089 Erfurt \\
Germany}
\email{dmitrii.legatiuk@uni-erfurt.de}
\subjclass{Primary 30G35, 60H40; Secondary 46G99, 60H15}

\keywords{Martingales, Brownian motion, It\^o formula, Stochastic Clifford analysis, Dirichlet problem}


\begin{abstract}
Clifford analysis has been the field of active research for several decades resulting in various methods to solve problems in pure and applied mathematics. However, the area of stochastic analysis has not been addressed in its full generality in the Clifford setting, since only a few contributions have been presented so far. Considering that the tools of stochastic analysis play an important role in the study of objects, such as positive definite functions, reproducing kernels and partial differential equations, it is important to develop tools for the study of these objects in the context of Clifford analysis. Therefore, in this work-in-progress paper, we present further steps towards stochastic Clifford analysis by studying random variables, martingales, Brownian motion, and It\^o formula in the Clifford setting, as well as their applications in Clifford analysis.
\end{abstract}

\maketitle

\section{Introduction}

Although stochastic differential equations play an essential role in real-world modelling processes, they have received only limited attention in the field of Clifford analysis. They play a vital role in real-world modelling processes scattering from diffusion and financial markets till describing phenomena from quantum mechanics, see \cite{da_Prato_1,Holden} for more details on modelling with stochastic differential equations. The approaches to the classical theory of stochastic equations comprise two general ways: (i) a semi-group approach utilising a semi-group generated by the differential operator of a stochastic DE, see e.g., \cite{da_Prato}, and (ii) Wick-product approach, where the classical objects of the stochastic calculus, such as e.g., Brownian motion, white noise, and It\^o integral, are transferred to the Wick setting, see \cite{Holden} for details. Independent on the way, analysis of stochastic DEs in a Clifford setting requires, first of all, a generalisation of the classical stochastic calculus to Clifford analysis.\par
Several authors have made first steps towards generalising tools of stochastic calculus to a hypercomplex setting in recent years. In \cite{Alpay_1}, positive definite functions have been studied in a quaternionic setting, where quaternionic random variables and stochastic processes have been introduced. Another group of works \cite{Alpay_2,Alpay_3,Alpay_4} has been devoted to studying stochastic calculus in general Grassmann algebras. Additionally, it is important to mention, that a connection between $\alpha$-hyperbolic harmonics and hyperbolic Brownian motion has been presented in \cite{Eriksson}. Nonetheless, a general Clifford setting has not been addressed so far.\par
Further ideas on generalising stochastic calculus to Clifford analysis have been presented in \cite{Bernstein}, where, among other results, Clifford random variables, Clifford white noise, and Clifford chaos expansion with the help of Hermite polynomials have been introduced. In this paper, we further extend and refine the results from \cite{Bernstein} by discussing martingales, Brownian motion, and It\^o formula in the Clifford setting. This way, we broaden ideas from \cite{nic,Getoor,Uboe} on complex Brownian motion and the It\^o formula to the hypercomplex setting. Additional motivation for generalising stochastic analysis to the Clifford setting comes from the fact, that tools of stochastic analysis provide a new point of view on studying classical deterministic problems and results in harmonic analysis, such as the solution of Dirichlet boundary value problems and properties of harmonic functions, see for example \cite{Oksendal,RogersWilliams}. In this paper we show that similar techniques can be used in the context of Clifford analysis by studying a Dirichlet problem for monogenic functions and proving the canonical Liouville's theorem with the help of the tools of stochastic Clifford analysis developed in this paper.\par
The paper is organised as follows: preliminaries from the classical Clifford analysis, as well as some supplementary results, are presented in Section~\ref{Section:Preliminaries}; Section~\ref{Section:Basic_stochastic_Clifford} presents basics of stochastic analysis in the Clifford setting, namely random variables, stochastic processes, and martingales; Section~\ref{Section:Brownian_motion} discusses Brownian motions in the Clifford setting, as well as shows its applications to some classical problems in Clifford analysis; finally, Section~\ref{Section:Ito_formula} is devoted to the It\^o formula and related results.\par

\section{Preliminaries and notations}\label{Section:Preliminaries}

Following \cite{Brackx,Guerlebeck_1}, we recall in this section some well-known facts from Clifford analysis. Let us consider the standard orthonormal basis $\{ \mathbf{e}_{1},\mathbf{e}_{2},\ldots,\mathbf{e}_{n} \}$ of the Euclidean vector space $\mathbb{R}^{n}$. From now on, $\mathcal{C}\ell(n)$ will denote the $2^{n}$-dimensional real Clifford algebra over $\mathbb{R}^{n}$ with the classical multiplication rules for basis vectors:
\begin{equation*}
\mathbf{e}_{j}^{2}=-1, \qquad \mathbf{e}_{j}\mathbf{e}_{k} + \mathbf{e}_{k}\mathbf{e}_{j} = 0, \quad j\neq k, \quad j,k=1,2,\ldots,n.
\end{equation*}
As usual, we identify $\mathbf{e}_{0}$ with the multiplicative identity $1$ of the real field. The Euclidean vector space $\mathbb{R}^{n+1}$ can be straightforwardly embedded in $\mathcal{C}\ell(n)$ by identifying the element $\mathbf{x}=(x_{0},x_{1},x_{2},\ldots,x_{n})$ with the Clifford para-vector $\mathbf{x}$ given by
\begin{equation*}
\mathbf{x} = x_{0} + \sum\limits_{k=1}^{n}x_{k}\mathbf{e}_{k} \in \mathcal{A}:=\mathrm{span}_{\mathbb{R}}\{1, \mathbf{e}_1, \ldots , \mathbf{e}_n\}.
\end{equation*}
Hence, the vector space  $\mathcal{A}$  is algebraically isomorphic to $\mathbb{R}^{n+1}$.\par
For an arbitrary para-vector $\mathbf{x}$, $\mathbf{Sc}(\mathbf{x}):= x_0$ and $\mathbf{Vec}(\mathbf{x}):= \sum_{k=1}^{n} x_k \mathbf{e}_k$ denotes its scalar part and its vector part, respectively. A basis of vector space for $\mathcal{C}\ell(n)$ is given by 
\begin{equation*}
\left\{ \mathbf{e}_0=1, \mathbf{e}_A = \mathbf{e}_{i_1} \cdots \mathbf{e}_{i_k}:  A=\{ i_1, \cdots, i_k \}, 1\leq i_1 < \cdots < i_k \leq n \right\},
\end{equation*}
implying that each element of $a\in\mathcal{C}\ell(n)$ (a Clifford number) can be represented in the form
\begin{equation*}
a=\sum_{A}a_{A}\mathbf{e}_{A}, \mbox{ where } a_{A}\in\mathbb{R}.
\end{equation*}
Additionally, the \textit{conjugation} in $\mathcal{C}\ell(n)$ is defined as the involutory anti-automorphism $\bar{\cdot}: \mathcal{C}\ell(n) \rightarrow \mathcal{C}\ell(n)$ given by its action on the basis elements: 
\begin{equation*}
\overline{a b} = \bar b ~\bar a, \quad a, b \in \mathcal{C}\ell(n), \mbox{ and } \bar{\mathbf{e}}_{0} = \mathbf{e}_{0}, \quad\bar{\mathbf{e}}_{k} = - \mathbf{e}_{k}, ~k=1, \cdots, n.
\end{equation*}
Thus, we have that $ \mathbf{x} \bar{\mathbf{x}} = \bar{\mathbf{x}}\mathbf{x} = \sum_{k=0}^n x_k^2 := |\mathbf{x}|^2$, is the Euclidean norm of $\mathbf{x}$ when identified with $x \in \mathbb{R}^{n+1}.$\par
Let now $\Omega$ be an open subset of $\mathbb{R}^{n+1}$ with a piecewise smooth boundary. A $\mathcal{C}\ell(n)$-valued function is a mapping
\begin{equation*}
f\colon\Omega\mapsto\mathcal{C}\ell(n) \mbox{ with } f(\mathbf{x})=\sum_{A}f_A(\mathbf{x})\mathbf{e}_{A}, \quad \mathbf{x}\in\Omega,
\end{equation*}
where the coordinates $f_A$ are real-valued functions defined in $\Omega$, i.e., $f_A\colon\Omega\to\mathbb{R}$, for all $A$. Continuity, differentiability and integrability of $f$ are defined coordinate-wisely. In the special case of para-vector valued functions we will denote them as $f(\mathbf{x})=\sum_{k=0}^{n}f^{k}(\mathbf{x})\mathbf{e}_{k}$.\par 
\begin{defn}
For continuously real-differentiable functions $f\colon\Omega\subset\mathbb{R}^{n+1}\to\mathcal{C}\ell(n)$, which is denoted for simplicity by $f\in C^{1}(\Omega,\mathcal{C}\ell(n))$, the operator
\begin{equation*}
D_{\mathbf{x}}:=\sum\limits_{k=1}^{n}\mathbf{e}_{k}\partial_{x_{k}}
\end{equation*}
is called the Dirac operator.
\end{defn}\par
Additionally, we would like also to consider the generalised Cauchy-Riemann operator in $\mathbb{R}^{n+1}, n\geq 1$, which is given by
\begin{equation*}
D = \partial _{x_0}+ \mathbf{e}_1 \partial_{x_1} + \ldots +\mathbf{e}_n  \partial_{x_n} = \partial_{x_0} + D_{\mathbf{x}},
\end{equation*}
and its conjugated operator
\begin{equation*}
\overline{D} = \partial_{x_0} - \mathbf{e}_1 \partial_{x_1} - \ldots -\mathbf{e}_n  \partial_{x_n} = \overline{ \partial_{x_0} + D} = \partial_{x_0} - D_{\mathbf{x}}.
\end{equation*}
Following \cite{Guerlebeck_Malonek_1999,Malonek_1990} and introducing the hypercomplex variables
\begin{equation*}
z_k = x_k - x_0e_k = -\frac{\mathbf{x}\mathbf{e}_k + \mathbf{e}_k\mathbf{x}}{2},
\end{equation*}
we can consider $n$ copies $\mathbb{C}_k$ of $\mathbb{C}$ by identifying $i\cong\mathbf{e}_k,$ $(k=1,\ldots,n), x_0 = \mathrm{Re} z, x_k \cong\mathrm{Im} z$, where $z\in \mathbb{C},$ and taking then $\mathbb{C}_k := -\mathbf{e}_k \mathbb{C}.$ Defining the Cartesian product $\mathbb{H}_k := \mathbb{C}_1 \times \cdots \mathbb{C}_k, (k=1,\ldots,n)$, we get the real linear vector space
\begin{equation*}
\mathbb{H}_n= \{ \mathbf{z} \,\colon \, \mathbf{z} = (z_1,\ldots, z_n) = (x_1,\ldots, x_n) - x_0(\mathbf{e}_1,\ldots , \mathbf{e}_n)\} \cong \mathbb{R}^{n+1} \cong \mathcal{A}.
\end{equation*}
The hypercomplex variables $z_{k}$ and the real linear vector space $\mathbb{H}_{n}$ introduced above will be used later in the paper for proving the It\^o formula in Clifford analysis.\par
\begin{defn}
A function $f\in C^{1}(\Omega,\mathcal{C}\ell(n))$ is called left (resp. right) monogenic in $\Omega$ if
\begin{equation*}
D f=0 \quad \mbox{in} \quad \Omega \quad (\mbox{resp.}, f D=0 \quad \mbox{in} \quad \Omega).
\end{equation*}
\end{defn}\par
For discussing basics of stochastic analysis in the Clifford setting, we need to introduce at first a proper Hilbert space structure, which is provided by a {\itshape right Hilbert-module} \cite{Brackx}:
\begin{defn}
A right (unitary) module over $\mathcal{C}\ell(n)$ is a vector space $V$ together with an algebra morphism $R\colon \mathcal{C}\ell(n) \rightarrow V$ (also called right multiplication), such that 
\begin{equation*}
R(ab+c)= R(b)R(a)+R(c).
\end{equation*}
In the particular case where $V=H\otimes \mathcal{C}\ell(n),$ with $\mathcal{H}$ a Hilbert space, we say that $V$ is a right Hilbert-module over $\mathcal{C}\ell(n)$.
\end{defn}
The inner product $(\cdot, \cdot)_{\mathcal{H}}$ in $\mathcal{H}$ gives rise to a {\itshape Clifford algebra-valued inner product} in $V$: 
\begin{equation*}
\langle f, g \rangle_{\mathcal{H}} := \sum_{A, B} (f_A, g_B)_{\mathcal{H}} \overline{\mathbf{e}}_{A} \mathbf{e}_{B},
\end{equation*} 
which is, strongly speaking, not the classical inner product due to lack of the positiveness axiom. Nonetheless, by restricting the Clifford algebra-valued inner product to its scalar part, the inner product can be defined.\par
By considering $\mu$-measurable Clifford algebra-valued functions in $\mathbb{R}^{n+1}$, the $L^{2}$-inner product can be stated as follows
\begin{equation*}
\langle f,g\rangle_{L^{2}(\mathbb{R}^{n+1},\mathcal{C}\ell(n))} = \int\limits_{\mathbb{R}^{n+1}} \overline{f(\mathbf{x})}g(\mathbf{x}) d\mu.
\end{equation*}\par

\section{Random variables, stochastic processes and martingales in the Clifford setting}\label{Section:Basic_stochastic_Clifford}

In this section we summarise some basic facts regarding stochastic analysis in the Clifford setting. Especially, a probability space, random variables, and stochastic process will be introduced. Some of the results presented in this section are built upon extending the work \cite{Bernstein}, where first steps towards stochastic Clifford analysis have been made. The final goal of this section is to introduce martingales in the Clifford setting. Moreover, for a better readability of the paper, we will provide all basic definitions needed for construction of a stochastic analysis in the Clifford setting, although some of these definitions are straightforwardly extended from the classical real case.\par
Similar to \cite{Alpay_1}, we start by defining the probability space and Clifford random variables:
\begin{defn}
The triple $(\Omega,\mathcal{F},P)$ is called a {\itshape probability space}, where $\Omega\subset\mathbb{R}^{n+1}$ is a given set, $\mathcal{F}$ is a $\sigma$-algebra on $\Omega$, and $P$ is a standard probability measure. A {\itshape Clifford random variable} $X$ is a $\mathcal{F}$-measurable function $X\colon \Omega\to \mathcal{C}\ell(n)$.
\end{defn}
It is important to remark, that since a Clifford random variable $X$ is a $\mathcal{C}\ell(n)$-valued function, then according to the discussion in Section~\ref{Section:Preliminaries}, this random variable can be expressed as follows
\begin{equation}
\label{Clifford_random_variable_expression_sum_real}
X=\sum_{A}X_{A}\mathbf{e}_{A},
\end{equation}
where $X_{A}$ are real-valued random variables. By help of the last representation, the classical proposition, see for example \cite{da_Prato_1}, can be straightforwardly generalised to the Clifford setting:
\begin{prop}
Let $X$ and $Y$ be Clifford random variables, then $\alpha X + \beta Y$ is a Clifford random variable for any $\alpha,\beta\in \mathcal{C}\ell(n)$.
\end{prop}\par
Next, following the standard theory of stochastic process, see for example \cite{Oksendal}, we call a probability measure $\mu_{X}$ on $\mathcal{C}\ell(n)$ induced by a random variable $X$
\begin{equation*}
\mu_{X}(\mathcal{F}) = P\left(X^{-1}(\mathcal{F})\right)
\end{equation*}
the {\itshape distribution of} $X$. Similarly, the number
\begin{equation*}
E[X] := \int\limits_{\Omega} X(\omega)dP(\Omega) = \int\limits_{\mathbb{R}^{n+1}} X d\mu_{X}(\mathbf{x}), \qquad (\mbox{if } \int\limits_{\Omega} |X(\omega)|dP(\Omega)<\infty)
\end{equation*}
is called {\itshape the expectation} of $X$ (w.r.t. $P$). By taking into account representation~(\ref{Clifford_random_variable_expression_sum_real}) of $X$, the expectation can also be written as follows
\begin{equation*}
E[X]=\sum_{A}E[X_{A}]\mathbf{e}_{A}.
\end{equation*}
Analogously to the classical theory \cite{da_Prato_1}, we define now a {\itshape stochastic Clifford process} as follows:
\begin{defn}
Let $V$ be a right Hilbert-module over $\mathcal{C}\ell(n)$, $(\Omega,\mathcal{F},P)$ be a probability space and let $I$ be an interval of $\mathbb{R}$. An arbitrary family $X=\left\{X(t)\right\}_{t\in I}$ of $V$-valued Clifford random variables $X(t)$, $t\in I$, defined on $\Omega$ is called a {\itshape stochastic Clifford process}. Additionally, we set $X_{t}(\omega)=X(t,\omega)$ for all $t\in I$ and $\omega\in\Omega$. Functions $X(\cdot,\omega)$ are called {\itshape trajectories} (or paths) of $X(t)$. Further, considering representation~(\ref{Clifford_random_variable_expression_sum_real}), a stochastic Clifford process can be written as follows
\begin{equation}
\label{Stochastic_process_Clifford_representation}
X(t) = \sum_{A}X_{A}(t)\mathbf{e}_{A}.
\end{equation}
\end{defn}
The parametrisation of Clifford random variables introduced in this definition is realised component-wisely for $\mathbf{x}\in\mathcal{C}\ell(n)$, similar to the case of vectors of random variables. Let us further introduce the following classical notions related to stochastic processes, which are straightforwardly generalised to the Clifford setting considering representation~(\ref{Stochastic_process_Clifford_representation}):
\begin{itemize}
\item A stochastic Clifford process $X(t)$ is said to be integrable (resp. square integrable) if
\begin{equation*}
E[|X_{A}(t)|]<\infty \, (\mbox{resp. } E[|X_{A}(t)|^{2}]<\infty) \mbox{ for all } t\in I \mbox{ and } A;
\end{equation*}
\item A stochastic Clifford process $X(t)$ is said to be bounded in $L^{p}$ if
\begin{equation*}
\sup\limits_{t\in I} E[|X_{A}(t)|^{p}]<\infty;
\end{equation*}
\item A stochastic Clifford process $X(t)$ is continuous if its trajectories are continuous.
\end{itemize}
Further definitions of regularity in a stochastic sense can also be adapted directly from the classical case.\par
Before introducing Clifford martingales, we need to recall the notion of a {\itshape filtration}:
\begin{defn}
\begin{itemize}
\item[(i)] A {\itshape filtration} on a probability space $(\Omega,\mathcal{F},P)$ is a family $\left\{\mathcal{F}_{t}\colon t\geq 0\right\}$ of $\sigma$-algebras such that $\mathcal{F}_{s}\subset \mathcal{F}_{t}\subset \mathcal{F}$ for all $s<t$.
\item[(ii)] A probability space together with a filtration is called a {\itshape filtered probability space}.
\item[(iii)] A stochastic Clifford process $\left\{X(t)\colon t\geq 0\right)\}$ defined on a filtered probability space with filtration $\left\{\mathcal{F}_{t}\colon t\geq 0\right\}$ is called {\itshape adapted} if $X(t)$ is  $\mathcal{F}_{t}$-measurable for any $t\geq 0$.
\end{itemize}
\end{defn}\par
Finally, we can introduce Clifford martingales:
\begin{defn}
Let $V$ be a right Hilbert-module over $\mathcal{C}\ell(n)$. An integrable $V$-valued stochastic Clifford process $\left\{X(t)\colon t\geq 0\right\}$ is a {\itshape martingale} with respect to a filtration $\left\{\mathcal{F}_{t}\colon t\geq 0\right\}$ if it is adapted to the filtration and
\begin{equation*}
E(X(t)|\mathcal{F}_{s}) = X(s),
\end{equation*}
almost surely for any pair of times $0\leq s\leq t$. The process is called a {\itshape submartingale} if $\geq$ holds, and a {\itshape supermartingale} if $\leq$ holds in the formula above.
\end{defn}\par
Next, for a fixed number $T>0$, let us denote by $\mathcal{M}_{T}^{2}(V)$ the space of all $V$-valued continuous, square integrable martingales $M$, such that $M(0)=0$. We have now the following proposition:
\begin{prop}
The space $\mathcal{M}_{T}^{2}(V)$, equipped with the norm
\begin{equation}
\label{Norm_martingale_space}
\|M\|_{\mathcal{M}_{T}^{2}(V)} := \left(E\left[\sup\limits_{t\in[0,T]}\|M(t)\|^{2}\right]\right)^{\frac{1}{2}},
\end{equation}
is a right-Banach module.
\end{prop} 
\begin{proof}
The proof of this proposition follows the standard arguments from the classical case, see for example \cite{da_Prato_1}. Therefore, we will only recall few basic ideas of the proof. At first, since we are in the Clifford setting, it is important to underline that the norm in~(\ref{Norm_martingale_space}) must be kept real-valued to ensure that $\|M(t)\|$ is a submartingale. Practically this implies, that in the case for example of inner product-induced norm, i.e. of a Hilbert space, the norm is constructed by using only the scalar component of the inner product. Thus, if $\|M(t)\|$ is a submartingale, then identity~(\ref{Norm_martingale_space}) defines a norm.\par
For proving completeness, we consider a Cauchy sequence $\left\{M_{n}\right\}$, which in stochastic setting implies
\begin{equation*}
E\left(\sup\limits_{t\in[0,T]} \|M_{n}(t)-M_{m}(t)\|^{2}\right) \to 0, \mbox{ as } n,m\to\infty.
\end{equation*}
After that, by using properties of Cauchy sequences and the Borel-Cantelli lemma, it is possible to show that $M$ is a continuous stochastic Clifford process. The final step of the proof is to recall that for a subsequence $\left\{M_{n_{k}}\right\}$ holds
\begin{equation*}
E(M_{n_{k}}(t)|\mathcal{F}_{s}) = M_{n_{k}}(s), 
\end{equation*}
almost surely if $0\leq s\leq t \leq T$ and $k\in\mathbb{N}$. After that, by letting $k$ tend to infinity, we get $E(M(t)|\mathcal{F}_{s}) = M(s)$ almost surely, implying that $M\in \mathcal{M}_{T}^{2}(V)$ and $M_{n}\to M$ in $\mathcal{M}_{T}^{2}(V)$.
\end{proof}\par
In Section~\ref{Section:Ito_formula}, a generalisation of the stochastic integral and It\^o formula to the Clifford setting will be discussed. To prepare this discussion, we need to introduce few more terms related to martingales, namely {\itshape local martingales} and {\itshape continuous semimartingale} \cite{RogersWilliams}:
\begin{defn}
A stochastic Clifford process $\left\{X(t)\colon t\geq 0\right\}$ is called a {\itshape local martingale}, if
\begin{itemize}
\item[(i)] $X(0)$ is $\mathcal{F}_{0}$-measurable,
\item[(ii)] $\left\{X(t)-X(0) \colon t\geq 0\right\}\in \mathcal{M}_{0,loc}$,
\end{itemize}
where $\mathcal{M}_{0,loc}$ is the space of local martingales null at $t=0$, i.e
\begin{equation*}
\mathcal{M}_{0,loc} := \{M \colon M(t), t\geq 0 \text{ is a continuous local martingale and } M_0= 0\}.
\end{equation*}
\end{defn}
Additionally to the space $\mathcal{M}_{0,loc}$, the following space needs to be introduced
\begin{equation*}
\mathbb{L}^2_{loc}(M) := \left\{ F \colon \begin{array}{l}
\left\{ F(t)\colon t\geq 0\right\} \text{ is progressively measurable } \\
\text{and } \int_0^{\infty} F^2 d\langle M \rangle < \infty
\end{array}  \right\},
\end{equation*}
where the term {\itshape progressively measurable} implies that a stochastic process is measurable with respect to the $\sigma$-algebra $\mathcal{B}([0,t])\otimes \mathcal{F}_{t}$ with $\mathcal{B}([0,t])$ denoting the Borel $\sigma$-algebra on $[0,t]$.\par
Finally, we need the notion of a {\itshape continuous semimartingale}:
\begin{defn}
A stochastic Clifford process $X(t), t\geq 0$ is a {\itshape continuous semimartingale} if it can be written as the sum
\begin{equation*}
X(t) = X(0) + M(t) + A(t),
\end{equation*}
where $X(0)$ is $\mathcal{F}_0$-measurable, $M$ is a continuous local martingale with $M(0)=0$ and $A(t)$ is a continuous adapted process with paths of locally finite variation with $A(0)=0$. The processes $M$ and $A$ are known as the martingale and finite variational parts of $X$, respectively.
\end{defn}\par
In the next section, after introducing Clifford Brownian motion, we will provide several examples of Clifford martingales and linking them to the Clifford Brownian motion.\par
\begin{rem}
As a summary of this section, we would like to underline that many constructions from the classical stochastic analysis can be directly transferred to the Clifford setting. The technical aspects of working with Clifford structures are first of all related to definitions of spaces, norms, and inner products, which are generally \textquotedblleft hidden\textquotedblright\, behind the classically formulated definitions.
\end{rem}\par

\section{Brownian motion and monogenic functions}\label{Section:Brownian_motion}

The aim of this section is to discuss Brownian motion in the Clifford context and relate it to monogenic functions. From now on we will only consider {\itshape white noise measure} $\mu_{X}$, i.e. the measure satisfying the Bochner-Minlos theorem, see \cite{Holden} for details. By the help of white noise measure, a {\itshape Clifford Brownian motion} can be introduced:
\begin{defn}\label{Definition:Clifford_Brownian_motion}
A para-vector-valued stochastic process $\left\{\mathbf{B}(t)\colon t\geq 0 \right\}$ of the form
\begin{equation*}
\mathbf{B}(t) := B_{0}(t) + \mathbf{e}_{1}B_{1}(t) + \ldots + \mathbf{e}_{n}B_{n}(t),
\end{equation*}
is called a {\itshape (linear) Clifford Brownian motion} with start $\mathbf{x}\in\mathbb{R}^{n+1}$, and where $B_{i}(t)$, $i=0,1,\ldots,n$ are classical one-dimensional Brownian motions, meaning that
\begin{itemize}
\item $B_{i}(0)=\mathbf{x}$;
\item the process has independent increments, i.e. for all times $0\leq t_{1}\leq t_{2}\leq \ldots \leq t_{k}$ the increments $B_{i}(t_{k})-B_{i}(t_{k-1})$, $B_{i}(t_{k-1})-B_{i}(t_{k-2})$, $\ldots$, $B_{i}(t_{2})-B_{i}(t_{1})$ are independent random variables;
\item for all $t\geq 0$ and $\sigma>0$, the increments $B_{i}(t+\sigma)-B_{i}(t)$ are normally distributed with expectation zero and variance $\sigma$;
\item almost surely, the function $t\mapsto B_{i}(t)$ is continuous.
\end{itemize}
Further, analogous to the classical case, we say that $\left\{\mathbf{B}(t)\colon t\geq 0 \right\}$ is a {\itshape standard Clifford Brownian motion} if $\mathbf{x}=0$.
\end{defn}
It is important to underline that only {\itshape para-vector-valued} Clifford Brownian motions are considered from now on. This restriction is necessary for providing clear connections between the Clifford Brownian motion and monogenic functions, which is presented in Theorem~\ref{thmDP}.\par
Next, we formulate the following theorem:
\begin{thm} 
Let $\{\mathbf{B}(t)\colon t\geq 0\}$ be a Clifford Brownian motion started at $\mathbf{x}\in\mathbb{R}^{n+1}$. Then the process $\{\mathbf{B}(t+s) - \mathbf{B}(s) \colon t,s >0\}$ is a Brownian motion started at the origin and is independent of $\{\mathbf{B}(t): 0\leq t \leq s \}$.
\end{thm}
\begin{proof} 
This is an immediate consequence of the independence of increments of Brownian motion. 
\end{proof}\par
In the sequel, we will use a Clifford Brownian motion adapted to a filtration. Let us briefly illustrate how such a Brownian motion can be constructed. Suppose $\{\mathbf{B}(t), t\geq 0\}$ is a Clifford Brownian motion defined on some probability space, then a filtration $\left\{\mathcal{F}_{t}^{0}\colon t\geq 0\right\}$ can be defined as follows
\begin{equation*}
\mathcal{F}_{t}^{0} = \sigma\left\{\mathbf{B}(s)\colon 0\leq s\leq t\right\},
\end{equation*}
which is the $\sigma$-algebra generated by the random variables $\mathbf{B}(s)$ for $0\leq s\leq t$. An important property of a Brownian motion is its independence on a filtration:
\begin{thm}\label{Theorem:Brownian_filtration_independent}
For all $s\geq 0$ the random process $\{\mathbf{B}(t+s) - \mathbf{B}(s), t\geq 0\}$ is independent of $\mathcal{F}_{t^{+}}$.
\end{thm}
\begin{proof} 
Because of the continuity of stochastic process $\{\mathbf{B}(t+s) - \mathbf{B}(s), t\geq 0\}$, the following equality holds for a strictly decreasing sequence $\{s_n\colon n\in \mathbb{N}\}$ converging to $s$:
\begin{equation*}
\mathbf{B}(t+s) - \mathbf{B}(s) = \lim_{n\to \infty} \mathbf{B}(s_n+t) - \mathbf{B}(s_n).
\end{equation*}
Then the Markov property implies that the right side of the above equation is independent of $\mathcal{F}_{t^{+}}$.
\end{proof}
Further we list two theorems, whose proofs are straightforward \cite{MP}, and, therefore, omitted:
\begin{thm}[Strong Markov property] 
For every almost surely finite stopping time $T$, the process $\{\mathbf{B}(T+t) - \mathbf{B}(T)\colon t\geq 0\}$ is a standard Clifford Brownian motion independent of $\mathcal{F}_{t^{+}}$.
\end{thm}
\begin{thm}[Reflection principle] 
If $T$ is a stopping time and $\{\mathbf{B}(t)\colon t\geq 0 \}$ is as standard Clifford Brownian motions, then the random process
\begin{equation*}
\mathbf{B}^*(t):= \mathbf{B}(t)\mathbb{1}_{\{t\leq T\}} + (2\mathbf{B}(T)-\mathbf{B}(t))\mathbb{1}_{\{t >T\}} =\left\{ \begin{array}{cl} \mathbf{B}(T) - \mathbf{B}(t), & \text{ if } t \leq T \\ \mathbf{B}(t), & \text{ if } 0\leq t \leq T \end{array} \right.
\end{equation*}
is a standard Clifford Brownian motions, and where $\mathbb{1}$ is the classical characteristic function of a set.
\end{thm}\par

\subsection{Clifford Brownian motion and Clifford martingales}

In this subsection, we will present two basic examples of Clifford martingales linking them to the Clifford Brownian motion introduced in Definition~\ref{Definition:Clifford_Brownian_motion}. These examples serve as a basis for future studies of the Clifford Brownian motion. We start with the following example:
\begin{ex}
For a Clifford Brownian motion $\left\{\mathbf{B}(t)\colon t\geq 0 \right\}$ we have
\begin{equation*}
E\left[\mathbf{B}(t)|\mathcal{F}_{s}^{+}\right] = E\left[\mathbf{B}(t)-\mathbf{B}(s)|\mathcal{F}_{s}^{+}\right] + \mathbf{B}(s) = E\left[\mathbf{B}(t)-\mathbf{B}(s)\right] + \mathbf{B}(s) = \mathbf{B}(s)
\end{equation*}
for $0\leq s\leq t$, where Theorem~\ref{Theorem:Brownian_filtration_independent} has been used. Hence, Clifford Brownian motion is a Clifford martingale, as expected.
\end{ex}
As the next example, we present the following lemma:
\begin{lem}
Let $\{\mathbf{B}(t)\colon t\geq 0\}$ be a Clifford Brownian motion adapted to a filtration $\mathcal{F}_{t^{+}}$, then the process
\begin{equation*}
\{\mathbf{B}^{2}(t) - t \colon t\geq 0\}
\end{equation*}
is a martingale.
\end{lem}
\begin{proof}
The proof is done by straightforward calculations and by using the fact, that the increments of Brownian motions are independent of a filtration. So, we can compute
\begin{equation*}
\begin{array}{ll}
& \displaystyle E\left[\mathbf{B}^{2}(t) - t |\mathcal{F}_{s}^{+} \right] = E\left[\mathbf{B}^{2}(t) - t |\mathcal{F}_{s}^{+} \right] \\
= & \displaystyle E\left[\left(\mathbf{B}(t)-\mathbf{B}(s)+\mathbf{B}(s)\right)^{2} - t |\mathcal{F}_{s}^{+} \right] \\
= & \displaystyle E\left[\left(\mathbf{B}(t)-\mathbf{B}(s)\right)^{2}|\mathcal{F}_{s}^{+} \right] + 2E\left[\mathbf{B}(t)-\mathbf{B}(s)|\mathcal{F}_{s}^{+} \right]\mathbf{B}(s) + E\left[\mathbf{B}^{2}(s)|\mathcal{F}_{s}^{+} \right] - t \\
& = \displaystyle t-s + \mathbf{B}^{2}(s) = \mathbf{B}^{2}(s) - s,
\end{array}
\end{equation*}
and thus, the lemma is proved. It is only important to underline, that since para-vector-valued Brownian motions are considered, the left and right multiplications must be carried through.
\end{proof}
The result of this lemma provides a basis for proving the second Wald's lemma, which goes beyond the scope of the current paper and will be addressed in future work.\par

\subsection{Clifford Brownian motion and monogenic functions}

Let now $B_{r}(\mathbf{x})$ denotes a ball of radius $r$ centred at $\mathbf{x}$, and $\partial B_{r}(\mathbf{x})$ be its boundary, as usual. Further, we will need the classical Hardy spaces $H^{p}\left(\partial \Omega\right)$ of $L^{p}$-functions on $\partial \Omega$ with $0<p<\infty$ and monogenic in $\Omega$. We start this section with the following theorem, which provides a connection between a Clifford Brownian motion and monogenic functions:
\begin{thm}\label{thmDP}
Let $\Omega$ be a domain, and let $(\Omega,\mathcal{F},P)$ be the probability space on which the Clifford Brownian motion $\{\mathbf{B}(t), t\geq 0\}$ started inside $\Omega$ is defined and suppose that $\left(\mathcal{F}_{t}\colon t\geq 0\right)$ is a filtration to which Clifford Brownian motion is adapted such that the strong Markov property holds. Further, let $\tau = \tau(\partial\Omega) = \min\left\{t\geq 0 \colon \mathbf{B}(t)\in \partial\Omega\right\}$ be the first hitting time of the boundary $\partial\Omega$. Let $\varphi\colon \partial\Omega\to\mathcal{C}\ell(n)$ be measurable and belong to $H^{p}\left(\partial \Omega\right)$, and such that the function $u\colon \Omega\to \mathcal{C}\ell(n)$ with
\begin{equation*}
u(\mathbf{x}) = E\left[\varphi\left(\mathbf{B}(\tau)\right)\mathbb{1}_{\tau<\infty }\right], \quad \mbox{for every } \mathbf{x}\in\Omega,
\end{equation*}
is locally bounded. Then $u$ is a monogenic function.
\end{thm}
\begin{proof}
For a ball $B_{\delta}(\mathbf{x})\subset\Omega$ let $\tilde{\tau}=\inf \left\{t>0 \colon \mathbf{B}(t)\notin B_{\delta}(\mathbf{x})\right\}$, then the strong Markov property implies that
\begin{equation*}
\begin{array}{lcl}
\displaystyle u(\mathbf{x}) & = & \displaystyle E\left[E\left[\varphi(\mathbf{B}(\tau))\mathbb{1}_{\tau<\infty} | \mathcal{F}^{+}(\tilde{\tau})\right]\right] \\
& = & \displaystyle E\left[u(\mathbf{B}(\tilde{\tau}))\right] = \int\limits_{\partial B_{r}(\mathbf{x})} u(\mathbf{y})\omega_{\mathbf{x},\delta}(d\mathbf{y}),
\end{array}
\end{equation*}
where $\omega_{\mathbf{x},\delta}$ is the uniform distribution on the sphere $\partial B_{\delta}(\mathbf{x})$. Therefore, $u$ possesses the mean value property, and considering that it has $\varphi\in H^{p}(\partial\Omega)$ as a boundary value, it is monogenic on $\Omega$. See \cite{Bernstein_1} for a detailed discussion on the relations between Hardy space $ H^{p}(\partial\Omega)$ and monogenicity in $\Omega$. 
\end{proof}\par
Next, we want to illustrate the application of tools constructed in this section for studying some known problems from the classical Clifford analysis. We start with the following Dirichlet problem for monogenic functions:
\begin{defn}[Dirichlet problem] 
Let $\Omega$ be a domain in $\mathbb{R}^{n+1}$ and let $\partial\Omega$ be its boundary, and assume further $u_0\in H^{2}(\partial\Omega)\cap C(\partial\Omega)$. Find a function $u\in C(\Omega)$, such that $u$ is monogenic on $\Omega$ and satisfies the boundary condition $u(\mathbf{x})=u_0(\mathbf{x})$.
\end{defn}
For discussing solvability of the Dirichlet problem, we need at first to recall the {\itshape Poincar\'e cone condition}:
\begin{defn}[Poincar\'e cone]
Let $\Omega\subset\mathbb{R}^{n+1}$ be a domain. We say that $\Omega$ satisfies the Poincar\'e cone condition at $\mathbf{x}\in\partial\Omega$ if there exists a cone $V$ based at $\mathbf{x}$ with opening angle $\alpha>0$, and $h>0$ such that $V\cap B_{h}(\mathbf{x})\subset \Omega^{c}$.
\end{defn}
It is worth to recall that the Lipschitz domains also satisfy the cone condition, see for example \cite{Adams}, and hence, the subsequent results are applicable for a large class of problems.\par
As the final preparation step for proving solvability of the Dirichlet problem, we need the following lemma:
\begin{lem}[\cite{MP}]\label{lem} 
Let $0<\alpha < 2\pi $ and $C_\mathbf{0}(\alpha) \subset \mathbb{R}^{n+1}$ is a cone based at the origin with opening angle $\alpha,$ and 
\begin{equation*}
a = \sup_{\mathbf{x}\in B_{\frac{1}{2}}(\mathbf{0})} P\{\tau(\partial B_1(\mathbf{0})) < \tau(C_\mathbf{0}(\alpha))\}.
\end{equation*}
Then $a<1$ and, for any positive integer $k$ and $h>0,$ we have 
\begin{equation*}
P\{ \tau(\partial B_{h}(\mathbf{z})) < \tau(C_\mathbf{z}(\alpha))\} \leq a^k ,
\end{equation*}
for all $\mathbf{x},\mathbf{z}\in\mathbb{R}^{n+1}$ with $|\mathbf{x}-\mathbf{z}|< 2^{-k}h$, where $C_\mathbf{z}(\alpha)$ is a cone based at $\mathbf{z}$ with opening angle $\alpha$.
\end{lem}
We provide a short prove of this lemma for the convenience of the reader:
\begin{proof}
Obviously $a < 1$. If $\mathbf{x}\in B_{2^{-k}}(\mathbf{0})$ then by the strong Markov property
\begin{equation*}
\begin{array}{c}
\displaystyle P\{ \tau(\partial B_{h}(\mathbf{z})) < \tau(C_\mathbf{z}(\alpha))\}  \\
\displaystyle \leq \prod_{i=0}^{k-1} \sup_{\mathbf{x}\in B_{2^{-k+i}}(\mathbf{0})}  P\{ \tau(\partial B_{2^{-k+i+1}}(\mathbf{0})) < \tau(C_\mathbf{0}(\alpha))\}  = a^k .
\end{array}
\end{equation*}
Therefore, for any positive integer $k$ and $h>0$, we have by scaling 
\begin{equation*}
P\{ \tau(\partial B_{h}(\mathbf{z})) < \tau(C_\mathbf{z}(\alpha))\} \leq a^k,
\end{equation*}
for all $\mathbf{x}$ with $|\mathbf{x}-\mathbf{z}|<2^{-k}h$.
\end{proof}
Now, we can formulate the main theorem for the Dirichlet boundary value problem:
\begin{thm}[Dirichlet problem] 
Suppose $\Omega \subset \mathbb{R}^{n+1}$ is a bounded domain with cone property, and $u_{0}$ is a continuous function on $\partial \Omega$. Let $\tau(\partial\Omega) = \inf\{t>0 \colon \mathbf{B}(t) \in\partial \Omega \}$. Then the function 
\begin{equation*}
u(\mathbf{x}) = E\left[u_0(\mathbf{B}(\tau(\partial U)))\right], \quad \text{for } \mathbf{x}\in \overline{\Omega},
\end{equation*}
is the unique continuous monogenic function on $\Omega$ with $u(\mathbf{x}) = u_0(\mathbf{x})\in H^{2}(\partial\Omega)\cap C(\partial\Omega)$ for all $\mathbf{x}\in \partial \Omega$.
\end{thm}
\begin{proof} 
The uniqueness of the solution follows from the uniqueness of monogenic functions. Considering that the boundary values belong to the Hardy space $H^{2}(\partial\Omega)$ of inner monogenic function, it follows that the function $u$ is monogenic by Theorem~\ref{thmDP}. We only need to prove that the Poincar\'e cone condition implies the boundary condition. For a fixed $\mathbf{y}\in \partial \Omega$ there is a cone $C_\mathbf{y}(\alpha)$ based at $\mathbf{y}$ with angle $\alpha >0$ such that $C_\mathbf{y}(\alpha)\cap B_h(\mathbf{y})\subset \Omega^c.$ Using Lemma~\ref{lem}, for any positive integer $k$ and $h>0$, we have 
\begin{equation*}
P\{ \tau(\partial B_{h}(\mathbf{z})) < \tau(C_\mathbf{z}(\alpha))\} \leq a^k
\end{equation*}
for all $\mathbf{x}$ with $|\mathbf{x}-\mathbf{y}|< 2^{-k}h$. Given $\varepsilon >0$, there is a $0<\delta \leq h$ such that $|u_0(\mathbf{v})-u_0(\mathbf{y})| < \varepsilon $ for all $\mathbf{v}\in \partial \Omega$ with $|\mathbf{v}-\mathbf{y}|< \delta$. For all $\mathbf{x}\in \overline{\Omega}$ such that $|\mathbf{y}-\mathbf{x}|<2^{-k}\delta$,
\begin{equation}
\label{eq}
\begin{array}{lcl}
\displaystyle |u(\mathbf{x})-u(\mathbf{y})| & = & \displaystyle \left|E_xu_0(\mathbf{B}(\tau(\partial \Omega)))- u_0(\mathbf{y})\right| \\
& \leq & \displaystyle E\left|u_0(\mathbf{B}(\tau(\partial \Omega)))- u_0(\mathbf{y})\right|. 
\end{array}
\end{equation}
If the Brownian motion hits the cone $C_\mathbf{y}(\alpha)$, which is outside the domain $\Omega$, before the sphere $\partial B_{\delta}(\mathbf{y})$, then $|\mathbf{y}-\mathbf{B}(\tau(\partial \Omega))| < \delta$, and $u_0(\mathbf{B}(\tau(\partial \Omega)))$ is close to $u_0(\mathbf{y})$. Hence, (\ref{eq}) is bounded from above by
\begin{equation*}
\begin{array}{lcl}
\displaystyle 2 \|u_0\|_{\infty} P_\mathbf{x}\{\tau(\partial B_{\delta}(\mathbf{y})) & < & \displaystyle \tau(C_\mathbf{y}(\alpha))\} +\varepsilon P_\mathbf{x}\{\tau(\partial \Omega) \\
& < & \displaystyle \tau(\partial B_{\delta}(\mathbf{y}))\} \leq 2 \|u_0\|_{\infty} a^k + \varepsilon,
\end{array}
\end{equation*}
where $\|\cdot\|_{\infty}$ is the classical $L^{\infty}$-norm. This implies that $u$ is continuous on $\overline{\Omega}$. 
\end{proof}\par
\begin{rem} 
If the domain $\Omega$ fulfils the cone condition, a solution of the Dirichlet problem can be simulated by running many independent Brownian motions, starting in $\mathbf{x}\in \Omega$ until they hit the boundary $\partial \Omega$ and letting $f(\mathbf{x})$ be the average of the values of $\varphi $ on the hitting points.
\end{rem}\par
Next, we show how the classical Liouville's theorem can be proved by the help of stochastic Clifford analysis:
\begin{thm}[Liouville's theorem] 
Any bounded monogenic function in $\mathbb{R}^{n+1}$ is constant.
\end{thm}
\begin{proof} 
Let $f$ be a monogenic function such that $|f(\mathbf{x})|< M < \infty $ for all $\mathbf{x}\in \mathbb{R}^{n+1}$. Further, let $\mathbf{x},\mathbf{y}$ be two distinct points in $\mathbb{R}^{n+1}$, and $H$ the hyperplane so that the reflection in $H$ takes $\mathbf{x}$ to $\mathbf{y}$. Let $\{\mathbf{B}(t) \colon t\geq 0\}$ be a Clifford Brownian motion started at $\mathbf{x}$, and $\{\mathbf{B}^*(t) \colon t \geq \}$ its reflection in $H$. Let $\tau(H)=\min\{t \colon \mathbf{B}(t)\in H \}$ and because of 
\begin{equation}
\label{bm*}
\{ \mathbf{B}(t) \colon t\geq \tau(H)\} \stackrel{d}{=}  \{\mathbf{B}^*(t) \colon t\geq \tau(H)\} .
\end{equation}
The monogenicity of $f$ implies that $E\left[f(\mathbf{B})(t))\right] = f(\mathbf{x})$ and decomposing into $t<\tau(H)$ and $t\geq \tau(H)$ we get
\begin{equation*}
u(\mathbf{x}) = E\left[u(\mathbf{B}(t) \mathbb{1}_{\{t<\tau(H)\}}\right] +  E\left[u(\mathbf{B}(t) \mathbb{1}_{\{t\geq\tau(H)\}}\right].
\end{equation*}
Using~(\ref{bm*}) we obtain
\begin{equation*}
\begin{array}{lcl}
\displaystyle \left| u(\mathbf{x})-u(\mathbf{y})\right| & = & \displaystyle \left| E\left[u(\mathbf{B}(t))\mathbb{1}_{\{t< \tau(H)\}}\right] - E\left[u(\mathbf{B}^*(t))\mathbb{1}_{\{t< \tau(H)\}}\right] \right|  \\
  & \leq & \displaystyle 2 M P\{t< \tau(H)\} \to 0, \quad \text{ as } t\to \infty. 
\end{array}
\end{equation*}
Thus $u(\mathbf{x})=u(\mathbf{y}),$ and since $\mathbf{x}$ and $\mathbf{y}$ were chosen arbitrarily, $u$ must be constant.
\end{proof}

\section{Stochastic integral and It\^o formula}\label{Section:Ito_formula}

We start this section by defining the It\^o integral:
\begin{defn} 
Let $X=M+A$ be a semimartingale with the martingale part $M\in \mathcal{M}_{0, loc}$ and the finite variational part $A$. Further, let $F\in \mathbb{L}^2_{loc}(M)$, then the It\^o integral of $F$ with respect to $X$ is the stochastic process of the form
\begin{equation*}
\int_0^* F(s) dX(s) = \int_0^* F(s) dM(s) + \int_0^* F(s) dA(s),
\end{equation*}
where the first terms on the right-hand side is an {\itshape It\^o integral}, and the second term is a (pathwise) {\itshape Lebesgue-Stieltjes integral}.
\end{defn}\par
Before introducing the It\^o formula, it is necessary to mention that we will use the notation e.g. $X_{i}(\cdot)$ implying $\mathbf{e}_{i}$ component of $n$-dimensional stochastic process. The It\^o formula is presented in the following theorem:
\begin{thm}[It\^o's formula, \cite{nic, RogersWilliams}]\label{Theorem:Ito} 
Let $f: [0,\infty)\times \mathbb{R}^n \to \mathbb{R}$ and suppose that the derivatives $f_t$ and $f_{x_ix_j}$ exist and are continuous for all $1\leq i,j \leq n.$ For $i=1, \ldots, n,$ suppose that 
\begin{equation*}
X_{i}(\cdot) = X_{i}(0) + M_{i}(\cdot) + A_{i}(\cdot) 
\end{equation*}
is a continuous semimartingale with the martingale part $M_i$ and the finite variational part $A_i$. Then, the following relation holds for all $t$
\begin{equation*}
\begin{array}{ll}
\displaystyle f(t, X(t)) - f(0, X(0)) = & \displaystyle \int_0^t f_t(s, X(s)) ds  \\
 & \displaystyle + \sum_{i=1}^n \int_0^t f_{x_i}(s, X(s)) dM_i(s) \\
 & \displaystyle + \sum_{i=1}^n \int_0^t f_{x_i}(s, X(s)) dA_i(s) \\
 & \displaystyle + \frac{1}{2} \sum_{i,j=1}^n \int_0^t f_{x_{i}x_{j}}(s, X(s)) d\langle M_i, M_j \rangle_s.  
\end{array}
\end{equation*}
\end{thm}
\begin{rem}[\cite{nic, RogersWilliams}] 
Very often, the It\^o's formula is abbreviated as follows
\begin{equation*}
\begin{array}{lcl}
\displaystyle df(t,X(t)) & = & \displaystyle f_t(t,X(t))dt + \sum_{i=1}^n f_{x_i}(t, X(t))dX_i(t) \\
& & \displaystyle + \frac{1}{2} \sum_{i,j=1}^n f_{x_{i}x_{j}}(s, X(s)) dX_i(t)dX_j(t)
\end{array}
\end{equation*}
with the convention that $dX_i(t) = dM_i(t) + dA_i(t)$, where
\begin{equation*}
dA_i(t) dA_j(t) = dA_i(t) dM_j(t) = 0, \mbox{ and } dM_i(t)dM_j(t) = d\langle M_i, M_j \rangle_t .
\end{equation*}
\end{rem}\par
Next lemma provides a generalisation of the It\^o formula to the Clifford setting:
\begin{lem}\label{Ito_formula_Clifford}
Let $f\colon \mathcal{C}\ell(n) \to \mathcal{C}\ell(n)$ be twice continuously differentiable, and let $X(t)=X_0(t)+\mathbf{e}_{1}X_{1}(t)+\ldots + \mathbf{e}_{n}X_{n}(t)$ be a continuous Clifford martingale, then the It\^o formula in the Clifford setting is given by
\begin{equation*}
\begin{array}{c}
\displaystyle df(t,X(t)) = \displaystyle  f_t(t,X(t))dt  + \int_{0}^{t} dZ_0(t)(D f) + dZ_1(t) \partial_{x_1}f + \ldots + dZ_n(t) \partial_{x_n}f   \\
 \displaystyle + \frac{1}{2}\int_{0}^{t} \sum_{i=0}^n dZ_0(t)dX_i(t) D\left(\frac{\partial f}{\partial x_i}\right) + \sum_{j=0}^n dZ_j(t)dX_i(t) \frac{\partial^2 f}{\partial x_j \partial x_i },
\end{array}
\end{equation*}
where $dZ_0(t) = dX_0(t)$ and $dZ_j(t) = e_j dX_0(t) + dX_j (t).$
\end{lem}
\begin{proof}
Proof of this lemma is based on the classical It\^o formula, and additionally, for proving this lemma we need to calculate the differential $\mathcal{C}\ell(n)$-valued 1- and 2-forms. Using the isomorphism between $\mathbb{H}_n$ and $\mathbb{R}^{n+1}$ we consider the mapping $f: \mathbb{R}^{n+1} \mapsto C\ell (n)$ as a function of the form
\begin{equation*}
f: \mathbb{H}_n \mapsto C\ell (n),
\end{equation*}
and its differential at $\mathbf{z}\in \mathbb{H}_n$ is then given by an $\mathbb{R}$-linear map $df_{\mathbf{z}}: \mathbb{H}_n\mapsto C\ell (n)$. By identifying the tangent space at each point of $\mathbb{H}_n$ with $\mathbb{H}_n$ itself, the differential of $C\ell (n)$-valued 1-form can be written as follows
\begin{equation*}
df = \partial_{x_0}f dx_0 + \partial_{x_1}f dx_1 + \ldots + \partial_{x_n}f dx_n.
\end{equation*}
Applying $x_0 = z_0$, $x_k = z_0\mathbf{e}_k + z_k$ (resp. $x_0 = z_0$, $x_k = \mathbf{e}_kz_0 + z_k)$, $k=1,\ldots , n$, we get
\begin{equation*}
\begin{array}{lcl}
df & = & (f D)dz_0 + \partial_{x_1}f dz_1 + \ldots + \partial_{x_n}f dz_n \\
& = & dz_0(D f) + dz_1 \partial_{x_1}f + \ldots + dz_n \partial_{x_n}f,
\end{array}
\end{equation*}
where the $n+1$ basic hypercomplex 1-forms $dz_k$ are defined by $dz_0 = dx_0$ and $dz_k = -\mathbf{e}_k dx_0 + dx_k$ for $k=1,\ldots, n$, respectively. It is well known that $d^{2}=0$ by definition. Nonetheless, for introducing a Clifford It\^o formula, we need to perform explicit computations of the differential 2-form by help of $df$:
\begin{equation*}
\begin{array}{lcl}
d^2f = d(df) & = & \displaystyle \sum_{j=0}^n \frac{\partial^2}{\partial x_j\partial x_0} f dx_jdx_0 + \sum_{i=1}^n\sum_{j=0}^n  \frac{\partial^2}{\partial x_{i}\partial x_{j}} f dx_jdx_i \\
& = & \displaystyle dz_0dx_0 D\left( \frac{\partial f}{\partial x_0}\right) + \sum_{j=0}^n dz_j dx_0 \frac{\partial}{\partial x_j}\left(\frac{\partial f}{\partial x_0}\right) + \\
& & \displaystyle + \sum_{i=1}^n \left( dz_0dx_i D\left(\frac{\partial f}{\partial x_i}\right) + \sum_{j=0}^n dz_j dx_i \frac{\partial}{\partial x_j}\left(\frac{\partial f}{\partial x_i}\right)\right) \\
& = & \displaystyle \sum_{i=0}^n \left(dz_0dx_i D\left(\frac{\partial f}{\partial x_i}\right) + \sum_{j=0}^n dz_jdx_i \frac{\partial^2 f}{\partial x_j \partial x_i }\right).
\end{array}
\end{equation*}
Comparing these calculations with the classical It\^o formula, we obtain the statement of lemma.\par
\end{proof}\par
Lemma~\ref{Ito_formula_Clifford} presents a general form of the It\^o formula. However, more specific forms of the It\^o formula, for example when a stochastic process is a Brownian motion \cite{Emery,Oksendal}, are also useful in practical applications. For connecting the discussion in Section~\ref{Section:Brownian_motion} with the It\^o formula, let us now consider the case that $\mathbf{X}(t)$ is a para-vector-valued standard Brownian motion, such that $X_0 = 0$. The process $\mathbf{B}(t) = \mathbf{x}+\mathbf{X}(t)$ is called Brownian motion starting from $\mathbf{x}$ (in general, it is possible to take for $B_0$ any random variable independent of the process $\mathbf{X}$). Since $\langle X_i,X_j\rangle_t =\delta_{ij} t$, we can write It\^o's formula, applied to a function $f$, as follows
\begin{equation*}
\begin{array}{lcl}
\displaystyle f(\mathbf{B}(t)) &  = & \displaystyle f(\mathbf{B}(0)) + \int_0^t \nabla f(\mathbf{B}(s)) \cdot d\mathbf{B}(s) + \frac{1}{2} \int_0^t \Delta f(\mathbf{B}(s)) ds \\
 & = & \displaystyle f(\mathbf{B}_0) + \int_0^t  dZ_0(t)(D f) + dZ_1(t) \partial_{x_1}f + \ldots + dZ_n(t) \partial_{x_n}f  \\
& & \displaystyle + \frac{1}{2} \int_0^t \Delta f(\mathbf{B}(s)) ds,
\end{array}
\end{equation*}
where $dZ_0(t) = dB_0(t)$ and $dZ_j(t) = e_j dB_0(t) + dB_j(t)$. If $f$ is a monogenic function the formula reduces to 
\begin{equation*}
f(\mathbf{B}(t)) = f(\mathbf{B}(0)) + \int_0^t  dZ_1(t) \partial_{x_1}f + \ldots + dZ_n(t) \partial_{x_n}f.  
\end{equation*}
\begin{rem}
It is important to remark, that the It\^o formula presented in Lemma~\ref{Ito_formula_Clifford} is a first attempt to develop this tool in the context of Clifford analysis. Moreover, in contrast to the classical results on conformal martingales and complex Brownian motions, see again \cite{nic,Uboe}, the It\^o formula in Lemma~\ref{Ito_formula_Clifford} does not significantly simplify itself in the case if $f$ is a monogenic function, as one would expect. Therefore, further analysis and studies are needed in this direction.
\end{rem}
\begin{rem}
Additionally, we would like to remark that all constructions presented in this section would work also for a general Clifford-valued function $f$, and not only for a para-vector-valued function considered above. However, because the Clifford Brownian motion is para-vector-valued, and it serves as a variable in the It\^o formula, it seems to be more appropriate to stay with para-vector-functions at the moment. Alternatively, a more general notion of the Clifford Brownian motion can be introduced, but this goes beyond the scope of current paper.
\end{rem}

\section{Summary and outlook}

In this paper, we have further developed ideas towards stochastic Clifford analysis by discussing random variables, stochastic processes, martingales, Brownian motion, and the It\^o formula in the Clifford setting. As a practical application of tools of stochastic Clifford analysis developed in this paper, we have illustrated how classical results in Clifford analysis, such as solvability of a Dirichlet problem for monogenic functions and Liuville's theorem, can be proved from the stochastic point of view. Further, the It\^o formula introduced in this paper, is a basis for addressing stochastic partial differential equations in Clifford analysis, which is the scope of our future work.\par


\subsection*{Acknowledgement}
We thank sincerely the reviewers for providing excellent constructive comments and suggestions, which have helped improving the paper.

\end{document}